\documentclass[letterpaper, 10 pt, conference]{ieeeconf}

\IEEEoverridecommandlockouts
\usepackage{cite}

\usepackage{amsmath,amsthm,amsfonts,amssymb,dsfont}
\usepackage[usenames]{color}
\usepackage{caption,subcaption}
\usepackage{graphicx}
\usepackage{epstopdf}
\usepackage{bbm}
\usepackage{enumerate}
\usepackage{mathtools}

\newcommand{\one}{{\mathbf{1}}}
\newcommand{\est}{{\hfill $\star$}}

\newcommand{\tran}{^{\top}}

\newcommand{\diag}{\mbox {\rm diag}}

\newcommand{\beq}{\begin{equation}}
\newcommand{\eeq}{\end{equation}}
\newcommand{\beqn}{\begin{equation*}}
\newcommand{\eeqn}{\end{equation*}}
\newcommand{\bea}{\begin{eqnarray}}
\newcommand{\eea}{\end{eqnarray}}
\newcommand{\beas}{\begin{eqnarray*}}
\newcommand{\eeas}{\end{eqnarray*}}
\newcommand{\ba}{\begin{array}}
\newcommand{\ea}{\end{array}}
\newcommand{\bit}{\begin{itemize}}
\newcommand{\eit}{\end{itemize}}
\newcommand{\ben}{\begin{enumerate}}
\newcommand{\een}{\end{enumerate}}
\newcommand{\dss}{\displaystyle}

\newcommand{\ap}[1]{ ^{ {\mathrm{#1} } }}

\newcommand{\Real}[1]{ { {\mathbb R}^{#1} } }
\newcommand{\Realp}[1]{ { {\mathbb R}_+^{#1} } }

\newtheorem{proposition}{Proposition}%[section]
\newtheorem{lemma}{Lemma}%[section]

\newcommand{\calD}{{\mathcal D}}

\newcommand{\calT}{{\mathcal T}}

\newcommand{\ve}{\varepsilon}

\title{Control of  Dynamic Financial Networks\\ (The Extended Version)} %

\author{Giuseppe Calafiore, %~\IEEEmembership{Fellow,~IEEE},
Giulia Fracastoro, %~\IEEEmembership{Member,~IEEE}
and Anton V. Proskurnikov %~\IEEEmembership{Senior Member,~IEEE}%
\thanks{The authors are with the Department of Electronics and Telecommunications, Politecnico di Torino, Turin, Italy. E-mails:\texttt{\{giuseppe.calafiore, giulia.fracastoro, anton.proskurnikov\}@polito.it}.}}

\begin{document}
%\begin{frontmatter}

\maketitle

%\author[DET]{Giuseppe C. Calafiore}\ead{giuseppe.calafiore@polito.it},
%\author[DET]{Giulia Fracastoro}\ead{giulia.fracastoro@polito.it},
%\author[DET]{Anton V. Proskurnikov}\ead{anton.p.1982@ieee.org},

%\address[DET]{Department of Electronics and Telecommunications, Polytechnic of Turin, Turin, Italy}  % Please supply

\begin{abstract}
The current global financial system forms a highly interconnected network where a default in one of its nodes can propagate to many other nodes, causing a catastrophic avalanche effect.
In this paper we consider the problem of reducing the financial contagion by introducing some targeted interventions that can mitigate the cascaded failure effects. We consider a multi-step dynamic model of clearing payments and introduce an external control term that represents corrective cash injections made by a ruling authority. The proposed control model can be cast and efficiently solved as a linear program. We show via numerical examples that the proposed approach can significantly reduce the default propagation by applying small targeted cash injections.
\end{abstract}
%\begin{keyword}
%TODO
%\end{keyword}
%\end{frontmatter}

%\date{}
%\maketitle

\section{Introduction}
In this paper we consider the problem of mitigating the effects of financial contagion via targeted and optimized interventions.
Recent studies~\cite{gale2007financial,battiston2010structure,GlasYou:16,elliott2014financial} highlighted the fact that in the current highly interconnected financial system, where banks and other institutions are linked via a network of mutual liabilities, a financial shock in one or few nodes of the network may  hinder the possibility for these nodes to fulfill their  obligations towards other nodes, and therefore provoke default. In turn, the nodes directly connected to the nodes that experienced the initial shocks receive reduced or no payments from these latter nodes, so their cash balances may be affected to the point of making impossible the fulfillment of their liabilities, hence of provoking further defaults, and so on in a cascaded fashion. The described mechanism may spread over the network as a contagion, provoking a possibly disastrous sequence of avalanche failures and defaults.

In the mainstream approach to the study of  default spreading in financial networks, see, e.g.,  \cite{EisNoe:01,GlasYou:16},
the contagion develops instantaneously, and in the aftermath of the contagion the nodes agree to settle for a set of mutual payments called {\em clearing payments} that brings the network to a new equilibrium after the shock. However, the assumption that all payments are simultaneous is quite unrealistic. For this reason, recently some works \cite{sonin2017banks,chen2021financial,banerjee2018dynamic,capponi2015systemic,ferrara2019systemic,kusnetsov2019interbank} proposed time-dynamic extensions of this model.
In particular, in~\cite{AutPaper} we consider a multi-step setting, in which defaults at one stage do not freeze all financial operations. Instead, in case of defaulted nodes, the residual claims are carried over to the next period, and so on until the end of the considered horizon. We show in~\cite{AutPaper} that multi-stage clearing payments can be computed by solving recursively a sequence of LP problems, and that  the multi-step setting may mitigate the cascaded failure effects by allowing shocks to be absorbed over time.

In this paper, we start from the setup of the multi-step model developed in~\cite{AutPaper} and introduce in the model an external control term
representing corrective cash injections at nodes to be performed by a ruling authority. The rationale is that a ruling authority, perhaps public, may intervene with minimal and targeted cash injections at certain nodes in order to prevent catastrophic cascaded failure events. We show that such control problem can be cast and efficiently solved as a linear program, an we provide numerical evidence of the fact that small targeted interventions at selected nodes (i.e., selected by the control algorithm itself) may suffice to avoid disastrous system-wide failures whose costs may be much larger than the amount necessary to prevent them. The notion of external injections of cash to reduce the contagion has been already deeply investigated \cite{minca2014optimal,amini2015control,amini2017optimal, fukker2021optimal}, in particular in the works on systemic risk measures (see, e.g., \cite{feinstein2017measures,biagini2019unified}). However, in most of the cases these models consider a single-step setting where all payments are simultaneous. %The problem of optimal control in a financial network has been already studied in \cite{minca2014optimal,amini2015control,amini2017optimal,fukker2021optimal,barratt2020multi}.
%Most of these works assumes that all payments are simultaneous.
Instead, \cite{barratt2020multi} considers a multi-step setting as in the proposed control problem. However, \cite{barratt2020multi} does not consider the presence of an external control term as in the proposed model. In addition, differently from the proposed problem, \cite{barratt2020multi} also assumes that entities cannot pay other entities more than the cash they have on hand. %This means that cash cannot make multiple steps through the network at once.
%However, to the best of the authors' knowledge, this is the first work that proposes a model with an external control term in a time-dynamic setting.

The paper is structured as follows. In Sec.~\ref{sec:EisNoe} we introduce the Eisenberg-Noe single-period model of a networked
financial system. Sec.~\ref{sec:multi} presents the proposed multi-step dynamic extension with an external control term.
Sec.~\ref{sec:control} introduces the problem of controlling the financial network by optimal cash injection.
Sec. \ref{sec:examples} shows two numerical examples. Conclusions are drawn in Sec.~\ref{sec:conc}.

\section{The Eisenberg-Noe model}
\label{sec:EisNoe}
We start by describing the classical  Eisenberg-Noe model of a networked financial system.
Consider  $n$ financial nodes (banks) who are subject
to  mutual liabilities  $\bar p_{ij}\geq 0$,
where $\bar p_{ij}$  represents
the payment due from node $i$ to node $j$.
The interbank liabilities constitute the {\em liability matrix}
$\bar P	\in\Real{n\times n}$, such that $[\bar P]_{ij}= \bar p_{ij}$ for $i\neq j=1,\ldots,n$, and
$[\bar P]_{ii}= 0$ for $i=1,\ldots,n$.
Also, nodes may receive cash from external entities, which are not part of the network, and we denote by
 $c\in\Realp{n}$ the vector whose $i$th component $c_i\geq 0$ represents the total cash in-flow from the external entities to  node $i$.
 Following the approach in  \cite{EisNoe:01} we further assume
 that payments towards  external entities are made  to a fictitious node that owes no liability to the other nodes (the corresponding row of liability matrix $\bar P$ is zero).
In the Eisenberg-Noe model time plays no role; specifically, all settlements of liabilities are assumed to be executed simultaneously at the
end of a fixed time period.
In  normal situations, at the end of the considered period each node $i$ is able to pay its liabilities in full, which means that each node $i$ receives an
inflow of liquidity $ \bar\phi_i\ap{in}\doteq  c_i+ \sum\nolimits_{k\neq i} \bar p_{ki}$ and pays out its liabilities by a total amount of
$\bar p_i \doteq  \bar\phi_i\ap{out}\doteq   \sum\nolimits_{k\neq i} \bar p_{ik}$.
%Due to limited liability, the outflow cannot exceed the inflow, that is it must hold that
%$
%\bar\phi_i\ap{in} - \bar\phi_i\ap{out}\geq 0$, $ i=1,\ldots,n
%$.
A critical
 situation instead occurs when (due to, e.g., a drop in the external  in-flow $c_i$) some bank $i$ cannot fully pay its debt.
 In this situation, the actual payments to other banks  have to be less than their nominal due values $\bar p_{ij}$.
We denote by  $p_{ij}\in [0,\bar p_{ij}]$, $i\neq j=1,\ldots,n$,  the {\em actual} inter-bank payments executed at the end of the period, which we  collect in matrix $P\in\Real{n,n}$.
Under the actual payments,  the cash  inflows and  outflows at each node $i=1,\ldots,n$, are respectively
$\phi_i\ap{in}\doteq c_i+ \sum\nolimits_{k\neq i} p_{ki}$,
$p_i \doteq \phi_i\ap{out} \doteq  \sum\nolimits_{k\neq i} p_{ik}$.
The vectors of inflows and outflows are thus
\begin{gather}
\phi\ap{in}\doteq c+P\tran\one,%\label{eq:inflow}\\
 \quad\phi\ap{out} \doteq p \doteq   P\one,\label{eq:in-outflow}
\end{gather}
where $\one$ denotes a vector of ones.
The nonnegative balance condition requires that
 $
 w \doteq \phi\ap{in}- \phi\ap{out}
 \geq 0$.
%\label{eq:networth}
%\eeq
%
A matrix  of mutual payments $P$, with $0\leq P\leq \bar P$,
 is said to be {\em admissible} if
  $w\geq 0 $.
  If the nominal liabilities $\bar P$ are admissible then payments $P = \bar P$  are such that all mutual obligations are met while maintaining the net worth of each node nonnegative, and no default arises.
 If instead $\bar P$ is not  admissible, then
 some nodes are in {\em default}, and all nodes
  must agree
  on a different set of admissible payments $P$, which
are upper bounded by $\bar P$, since no node should pay more than due.
Moreover, when a node is in default, it must pay out all of its cash inflow to the creditor nodes: each node $i$ pays out $\bar p_i$ or  pays out its whole inflow $\phi_i\ap{in}$. Therefore, a {\em clearing  payment matrix} $0\leq P \leq \bar P$ obeys the relation
\beq
P\one = \min(\bar P\one , c+ P\tran\one).
\label{eq:clearingM}\eeq
One clearing matrix  $0\leq P \leq \bar P$ satisfying~\eqref{eq:clearingM} can be found~\cite{Journal2021} by solving an optimization problem of the form
\beq
\min_P\, f(P),\quad \mbox{s.t.: } 0\leq P \leq \bar P,\; P\one  \leq  c+ P\tran\one,
\label{eq:clearingoptP} \eeq
where $f(P)$ is any decreasing function of the matrix argument $P$ on $[0,\bar P]$, that is, a function such that $\bar P\geq P^{(2)} > P^{(1)}\geq 0$, $P^{(2)} \neq P^{(1)}$, implies
$f(P^{(2)}) < f( P^{(1)})$.
Possible choices for  $f$  in~\eqref{eq:clearingoptP} are for instance
$f(P) = \|\bar \phi\ap{in} -\phi\ap{in}\|_1$ and
$f(P) = \|\bar \phi\ap{in} -\phi\ap{in}\|_2^2$, where
 $\phi\ap{in}(P)  = c +P\tran \one$. The optimal solution of~\eqref{eq:clearingoptP}, however, is in general non unique.

In practice, payments under default are subject to further regulations. A commonly used ``local fairness'' rule is that
the outstanding claims should be redistributed
 based on a proportionality (pro-rata)  rule.
We define the relative proportion of payment due nominally by node $i$ to node $j$  as
 \beq
 a_{ij}\doteq \left\{\ba{ll}\dss \frac{\bar p_{ij}}{\bar p_i}  & \mbox{if }  \bar p_i> 0 \\
 1 & \mbox{if }  \bar p_i=0\,\mbox{ and } i=j \\
 0 & \mbox{otherwise}.
 \ea\right.
 \label{eq:prorataA}
 \eeq
 Computing these proportions for all $i,j$ we  form the {\em relative liability} matrix $A = [ a_{ij}]$.
     By definition,  matrix $A$ is row-stochastic, that is  $A\one = \one$.
The so called {\em pro-rata rule}   imposes that payments are due in proportion  to
the rates fixed in matrix $A$, that is
 $ p_{ij}= a_{ij}  p_i ,\quad \forall i,j$,  where $p_i$ is the out-flow defined in~\eqref{eq:in-outflow}.
Since $p\doteq P\one$, the pro-rata rule imposes a set of linear equality constraints
on the entries of $P$, namely $P = \diag(P\one)A =  \diag(p)A $.
Under the pro-rata rule, the full payment matrix $P$ is determined by vector $p$; problem~\eqref{eq:clearingoptP} simplifies in this case to
\beq
\min_P\, f(p),\quad \mbox{s.t.: } 0\leq p \leq \bar p,\; p  \leq  c+ A\tran p,
\label{eq:clearingoptP_pr} \eeq
and it holds that for any decreasing $f$ the solution $p^*$ to \eqref{eq:clearingoptP_pr} is unique and it represents a \emph{clearing vector}, that is, it satisfies $p = \min(\bar p , c+ A\tran p)$, see
e.g.~\cite[Lemma~1]{Journal2021}.

Even though most of the works on financial contagion impose the pro-rata rule, in this paper we consider also the more general case without such constraint. The non-proportional clearing mechanism can significantly reduce the impact of a financial shock \cite{Journal2021,barratt2020multi}. In addition, it may also be extended to promote virtuous behaviors such as rescue consortium \cite{rogers2013failure}.

\section{A multi-stage model with controls}
\label{sec:multi}
%We next extend the single-stage model described in the previous section to a multi-period dynamic setting.
 As already observed,  the  default and clearing model discussed in the previous section, which coincides with the mainstream one studied in the literature~\cite{GlasYou:16},
 is a single-period model, meaning  that the described process assumes that at one point in time (the end of a fixed period), all liabilities are claimed and due simultaneously, and that the entire network of banks becomes aware of the claims and possible defaults and instantaneously agrees on the clearing payments. All financial operations of  defaulted nodes are  frozen, which possibly induces propagation of the default to other neighboring nodes, in an avalanche fashion, see, e.g.,~\cite{Massai:2021}. In~\cite{AutPaper}, we propose a dynamic multi-step model in which financial operations are allowed for a given number of time periods after the initial theoretical defaults (here named {\em pseudo-defaults}). In this way  some nodes may actually {\em recover} and eventually manage to fulfill their obligations by the end of the allotted time horizon.  We next describe the dynamic model from~\cite{AutPaper}, and introduce into this model additional control inputs that were not considered in~\cite{AutPaper}.

We consider a discrete time horizon $t = 0,1,\ldots,T$, with periods of fixed duration (e.g., one day, or one month, etc.), where $T\geq 0$
represents the final time of the horizon. For brevity, denote $\calT\doteq\{0,\ldots,T-1\}$. Throughout the text, a sequence of vectors or matrices $(f(t),\,t\in\calT)$ is denoted by
\[
[f]\doteq (f(0),\ldots,f(T-1)).
\]

Extending to the multi-stage case  the basic model described in Section~\ref{sec:EisNoe}, we let
$\bar P(t)	\in\Real{n\times n}$,  and $ P(t)	\in\Real{n\times n}$, $t\in\calT$ denote the
nominal liabilities matrices and the actual payment matrices at time $t$, respectively.
We let $c(t) = e(t) + u(t) \geq 0$ denote the sum of the  vector $e(t)\geq 0$  of cash inflows at nodes from the external sector at time $t$, plus
the vector $u(t)\geq 0$
of additional ``control inflows'' injected at nodes
at time $t$ by the control authority.
We further define $\bar P \doteq \bar P(0)$ as the matrix of initial liabilities, and we let the pro-rata matrix $A$ be defined
as in \eqref{eq:prorataA} according to these initial liabilities.
In this work we can deal indifferently with models with full payment matrices, or with
matrices constrained by the pro-rata condition: in this latter case we shall simply include the linear equality constraint
$P(t) = \diag(p(t))A$ on the payment matrices, where $p(t) \doteq P(t)\one$.
 The net worth $w_i(t)$ of node $i$ at  $t$ evolves in accordance to
$
 w_i(t+1) = w_i(t) + \phi_i\ap{in} (t)  - \phi_i\ap{out}(t)$,
% \label{eq_wdynamics}
 %\]
 or, in the vector form,
 \beq\label{eq:wdynamics2free}
 w(t+1) = w(t) + c (t) + P\tran(t)\one - P(t)\one,\quad t\in\calT.
 \eeq
Similar to the single-period case discussed in Section~\ref{sec:EisNoe}, the limited liability condition requires that
$w(t)\geq 0$ at all $t$. It may therefore happen that  a payment $p_{ij}(t)<\bar p_{ij}(t)$ in order to guarantee  $w_i(t) \geq 0$.
When this happens at some $t< T$, instead of declaring default and freezing the financial system, we allow operations to continue up to the final time $T$, updating the
due payments according to the  equation
$
\bar p_{ij}(t+1) =  \alpha \left(\bar p_{ij} (t) - p_{ij} (t)  \right) $,
%  \label{eq_pdynamics}
 where $\alpha \geq 1$ is the interest rate applied on past due payments.
This can be written as
 \beq
\bar P(t+1) = \alpha \left(\bar P (t) - P (t)  \right), \quad t\in\calT.
  \label{eq_pdynamics}
 \eeq
 The recursions~\eqref{eq:wdynamics2free} and~\eqref{eq_pdynamics} are
  initialized with
 $
 w(0)  = 0$, $\bar P(0 ) =  \bar P$, where  $\bar P$ is the initial liability matrix.
  The meaning of equation~\eqref{eq_pdynamics} is that if  a due payment at $t$ is not paid in full, then the residual debt is added to the nominal liability for the next period,
 possibly increased by an interest factor $\alpha \geq 1$. This mechanism allows for a node which is technically in default at a time $t$ to continue operations and (possibly) repay its dues in subsequent periods.
 Notice that matrix $\bar P(t)$ is time-varying and depends on the actual payment matrices $P(0),\ldots,P(t-1)$; the final matrix $\bar P(T)$ contains the residual debts at the end of the final period.
 %
%As long as $w_i(t) > 0$ the dues of $i$ are payed in full, that is $ p_{ij} (t) = \bar p_{ij} (t)$,
% otherwise, the dues are repaid partially, that is $ p_{ij} (t) \leq \bar p_{ij} (t)$.
 %

 The payment matrices $P(t)$ are  subject to the constraints%, for $t\in\calT$
\bea
  && 0\leq  P(t) \leq \bar P(t)\quad\forall t\in\calT,  \label{eq.cond-p-1}\\
  && P(t) \one \leq   w(t) + c (t) + P(t)\tran \one\quad \forall t\in\calT, \label{eq.cond-p-2}
\eea
where \eqref{eq.cond-p-1} represents the requirement that actual payments never exceed the nominal liabilities, and
\eqref{eq.cond-p-2} represents the requirement that $w(t+1)$, as given in \eqref{eq:wdynamics2free}, remains nonnegative at all $t$.
%
%\color{blue}
Conditions  \eqref{eq.cond-p-1}, \eqref{eq.cond-p-2}  can be made explicit by eliminating the variables $w(t)$ and $\bar P(t)$, which
by using \eqref{eq:wdynamics2free}--\eqref{eq_pdynamics} can be expressed as
\begin{gather}
 \bar P(t) = \alpha^t \bar P(0) -  \sum\nolimits_{k=0}^{t-1} \alpha^{t-k}P(k), \label{eq:Pbart} \\
 w(t) = C(t-1) + \sum\nolimits_{k=0}^{t-1} \left(P\tran(k) - P(k) \right)\one , \label{eq.w-evolves} \\
 %\bar P(t)  - P(t)  = \alpha^t \bar P(0) -  \sum_{k=0}^{t} \alpha^{t-k}P(k), \; t\in [1:T-1]\notag,
%where
% \beq
 C(t) \doteq \sum\nolimits_{k=0}^{t} c(k)
 =  \sum\nolimits_{k=0}^{t} [e(k) + u(k)].
  \label{eq:cumulativeinflow}
  \end{gather}
 %\eeq
 %
Conditions~\eqref{eq.cond-p-1}, \eqref{eq.cond-p-2} can thus be rewritten as
\begin{gather}
P(t)\geq 0,\label{eq.cond-p-1+}\;\;\forall t\in\calT \\ % \notag\\
\sum\nolimits_{k=0}^{t} \alpha^{t-k}P(k) \leq \alpha^t \bar P\;\;\forall t\in\calT\label{eq.cond-p-1a} \\%\notag\\
C(t) + \sum\nolimits_{k=0}^{t} \left(P(k) \tran- P(k) \right)\one \geq 0\;\forall t\in\calT\label{eq.cond-p-2a}
\end{gather}
In the case when pro-rata is enforced
the above conditions can be rewritten in terms of
the outflow vectors only:
\begin{gather}
p(t)\geq 0,\label{eq.cond-p-1++}\\
\sum\nolimits_{k=0}^{t} \alpha^{t-k}p(k) \leq \alpha^t \bar p\label{eq.cond-p-1b}\\
C(t) + \sum\nolimits_{k=0}^{t} \left(A\tran p(k)- p(k) \right) \geq 0\label{eq.cond-p-2b}\\
\forall t\in\calT,\quad \bar p \doteq \bar P \one.\notag
\end{gather}
We say that a payment sequence $[P]$ is feasible for a {\em given} $[c]$, if it satisfies
\eqref{eq.cond-p-1+}--\eqref{eq.cond-p-2a}. Analogously, under pro-rata, a sequence $[p]$
is feasible for given  $[c]$, if it satisfies
\eqref{eq.cond-p-1++}--\eqref{eq.cond-p-2b}.

We now introduce a system-level cost criterion, based on the total difference between the nominal and actual payments at the nodes. We define the \emph{loss} at period $t$ by
\beq\label{eq.delta}
\begin{aligned}
\delta(t) &\doteq & \sum_{i,j=1}^n \left( \bar p_{ij}(t) - p_{ij}(t) \right)  =
 \one\tran \left( \bar P(t) - P(t) \right)  \one.
 \end{aligned}
\eeq
Observe that $\delta(t) \geq 0$ for all $t$ and for any feasible $[P]$, and
$\delta(t) = 0$ if and only if $P(t) = \bar P(t)$, that is when no default occurs.
Therefore, $\delta(t) $ can be taken as a measure of the effects of defaults at stage $t$.
The overall cost function is then defined as the \emph{total} loss over the time horizon
\beq
L([P]) \doteq  \sum_{t=0}^{T-1} \delta(t)  = a_0 \one\tran \bar P \one -  \sum_{t=0}^{T-1}a_t \one\tran P(t) \one,
\label{eq:Loss}
\eeq
where we derive the second equality from~\eqref{eq:Pbart} and define the constants $a_0>a_1>\ldots>a_{T-1}$ as
\beq
a_t\doteq \sum_{j=0}^{T-t-1} \alpha^j = \left\{ \ba{ll}
 \frac{\alpha^{T-t}-1}{\alpha -1}, & \mbox{if } \alpha > 1\\
T-t, & \mbox{if } \alpha = 1.
\ea \right.
\label{eq:ai}
\eeq
Now, for \emph{fixed} input flows $[c]$, the multi-stage clearing payments
$[P^*]$, are defined as the optimal solution to
\beq
\min_{[P]} L([P])\quad \mbox{s.t.:} \; \eqref{eq.cond-p-1+}-\eqref{eq.cond-p-2a}. %\eqref{eq.cond-p-1++}-\eqref{eq.cond-p-2b}.
\label{eq:prob01a}
\eeq
Analogously, under the pro-rata rule, the clearing payments are defined via vectors $p^*(t)$, $t\in\calT$ which solve the LP
\beq
\min_{[p]} L([p])\quad \mbox{s.t.:} \; \eqref{eq.cond-p-1++}-\eqref{eq.cond-p-2b},
\label{eq:prob01b}
\eeq
where $p(t) = P(t) \one$, and $L([p])\doteq %\sum_{t=0}^{T-1}\one\tran(\bar p(t) -p(t))$.
  a_0 \one\tran \bar p -  \sum_{t=0}^{T-1}a_t \one\tran p(t) $.
The properties of the optimization problems \eqref{eq:prob01a} and \eqref{eq:prob01b} (which are in fact LP) have been studied in~\cite{AutPaper}.

 \section{Multi-step control of the dynamic network}
\label{sec:control}
We next consider the problem of controlling the financial network by optimal injections of cash $[u] =(u(0),\ldots,u(T-1))$ at nodes.
The cumulative  amount of cash injected from $0$ to  $t$ is
\beq
B(t) \doteq \sum_{\tau=0}^{t} \one\tran u(\tau), \quad t\in\calT.\label{eq.budget}
\eeq
The control objective we propose to minimize is defined as
\beq\label{eq:Ja}
J([P],[u]) = (1-\eta)L([P]) + \eta\one\tran \bar P(T) \one +\gamma B(T-1),
\eeq
where $L([P])$ is given in \eqref{eq:Loss}, $\gamma \geq 0$ is a given penalty on
the total control cash, and $\eta\in[0,1]$ is a
 weight on the terminal cost ($\bar P(T)  = 0$ if and only if there is no default at the terminal time).
The control problem is then stated as
\bea
\min_{[P],[u]} & J([P],[u])  \label{eq:control1a}\\
\mbox{s.t.:} & [P] \geq 0,\quad [u]\geq 0 \nonumber \\
& \sum\nolimits_{k=0}^{t} \alpha^{t-k}P(k) \leq \alpha^t \bar P   ,\;   t\in\calT \label{eq.cond1a}\\
& C(t) + \sum\nolimits_{k=0}^{t}  \left(P(k) \tran- P(k) \right)\one \geq 0, \; t\in\calT  \label{eq.cond2a}\\
& B(t) \leq F(t), \;t\in\calT  \label{eq.cond3a}
\eea
where $C(t)$ is given by \eqref{eq:cumulativeinflow} for fixed $[e]$, and $F(t)\geq 0$ is a given \emph{nondecreasing
sequence} that represents the maximum budget available up to time $t$ for controlling the network.

Under the pro-rata rule  the control problem simplifies to
\bea
\min_{[p],[u]} & J([p],[u])  \label{eq:control1}\\
\mbox{s.t.:} & [p] \geq 0,\quad [u]\geq 0 \nonumber \\
& \sum\nolimits_{k=0}^{t} \alpha^{t-k}p(k) \leq \alpha^t \bar p   ,\;   t\in\calT \label{eq.cond1}\nonumber\\
& C(t) + \sum\nolimits_{k=0}^{t}  \left( A\tran p(k) - p(k) \right) \geq 0, \; t\in\calT  \label{eq.cond2}\nonumber\\
& B(t) \leq F(t), \;t\in\calT ,  \label{eq.cond3}
\eea
where the cost function is
\beq\label{eq:J}
J([p],[u]) = (1-\eta)L([p]) + \eta\one\tran \bar p(T)   +\gamma B(T-1),
\eeq

%\eeq

In the case where $\eta\in [0,1)$ (the loss accumulated over time is penalized) and $\gamma>0$ (the total control cash is penalized),
the solutions to the problems~\eqref{eq:control1a} and~\eqref{eq:control1} enjoy a number of important properties. Denote the optimal sequences of payment matrices and control inputs by $P^*(t)$ and $[u^*]$ respectively (in the problem~\eqref{eq:control1} $P^*(t) = \diag(p^*(t))A$). To each optimal solution, we associate the sequences $p^*(t)=P^*(t)\one$ and
$[{\bar P}^*],[{\bar p}^*], [c^*],[w^*],[B^*],[\delta^*]$.
\begin{lemma}\label{lem.1}
Let $\eta\in [0,1),\gamma>0$. Then, all optimal processes in problems~\eqref{eq:control1a} and~\eqref{eq:control1} enjoy
the following properties:
\begin{enumerate}
\item The absolute debt priority rule is respected:
\beq\label{eq.no-defferal}
p^*_i(t)=\min\left(\bar p_i^*(t),w^*(t)+c^*_i(t)+\sum\nolimits_{j\ne i}p_{ji}^*(t)\right)
\eeq
for all $i=1,\ldots,n$ and $t\in\calT$;
\item A bank utilizes the injected liquidity immediately by paying out all its balance:
if $u_i^*(t)>0$, then $w_{i}^*(t+1)=0$ and, moreover, $w_i(s)=0\,\forall s\leq t$.
\item If $B^*(t_*)<F(t_*)$ at some period $t_*<T-1$, then no liquidity is injected after period $t_*$: $u^*(t)=0\,\forall t>t_*$.
\end{enumerate}
\end{lemma}

Notice that the first property shows that the optimal clearing policy prohibits unnecessary deferrals of payments: bank $i$ pay out its liability $\bar p_i^*(t)$ as soon as this is possible.

Note also that if the whole control budget is available at $t=0$, i.e., $F(0)=\ldots=F(T-1)$, then we either have $t_*=0$ and $B^*(0)<F(0)$ or the whole budget is used at time $t=0$. In both situations, one obviously has $u^*(t)=0\,\forall t\geq 1$. Similarly, if $F(k)=\ldots=F(T-1)$, then there are no control actions after period $k$: $u^*(t)=0\,\forall t\geq k+1$.

\subsection{Dealing with uncertainty in the external payments}
In problem  \eqref{eq:control1}  we assumed that the whole stream $[e] = (e(0),\ldots,e(T-1))$ of cash inflows from the external sector to the nodes is precisely known in advance. In this section we consider instead a more realistic scenario in which the inflows are known only up to some interval of uncertainty. More precisely, we assume that
\[
e(t) = \hat e(t) + d(t), \quad t=0,\ldots,T-1,
\]
where $ \hat e(t)\geq 0$ is the nominal predicted value of the external cash inflow at $t$, and $d(t)$ is an unpredictable uncertainty on this value, assumed to bounded in magnitude so that
\[
|d_i(t)|\leq r_i(t) \doteq \epsilon(t) \hat e_i(t) , \quad i=1,\ldots,n;\; t=0,\ldots,T-1,
\]
 where $ \epsilon(t) \in(0,1)$ is the given relative error level at $t$.
 We let %$[d]$ denote the collection $(d(0),\ldots,d(T-1))$ and we let
 $\calD$ denote the uncertainty set on $[d]$, that is $\calD = \{[d] : |d_i(t)|\leq r_i(t) , i=1,\ldots,n; t=0,\ldots,T-1\}$.

 For simplicity of exposition and notation we treat here only the case of proportional payments, which simplifies the problem and allows us to deal only with vector variables $p(t)$ instead of matrix variables $P(t)$. The whole reasoning reported below, however, carries over to the matrix case with only formal and notational modifications.

 The decision variables $[p],[u]$ of  problem \eqref{eq:control1} are next assumed to be prescribed by a reactive policy that allows adjustments in consequence to deviation of the external inflows from their nominal values: for all $t\in\calT$ we let
 \begin{align}
 p(t) &=& \hat p(t) +\Theta(t) (e(t) - \hat e(t) )=\hat p(t) +\Theta(t)d(t)  \label{eq:policy}\\
 u(t) &=& \hat u(t) +\Gamma(t) (e(t) - \hat e(t) )=\hat u(t) +\Gamma(t)d(t)\notag.
 \end{align}
 and $[\hat p], [ \hat u]$ are now the new decision variables, together with the collections of reaction matrices $[\Theta], [\Gamma]$.
 The control problem \eqref{eq:control1} is now cast in a worst-case setting as follows
 \bea
\min_{[\hat p],[\hat u], [\Theta], [\Gamma]} & \max_{[d]\in\calD} J([p],[ u])  \label{eq:control1wc} \\
% J([\hat p],[\hat u], [\Theta], [\Gamma],[d])  \label{eq:control1wc}\\
\mbox{s.t.:} & \min_{[d]\in\calD} [p] \geq 0, \quad \min_{[d]\in\calD}[u]\geq 0 \nonumber \\
& \max_{[d]\in\calD} \sum\limits_{k=0}^{t} \alpha^{t-k}p(k) \leq \alpha^t \bar p,\,t\in\calT  \label{eq.cond1wc}\nonumber\\
& \min\limits_{[d]\in\calD}  C(t) + \sum\limits_{k=0}^{t} \left(A\tran p(k)- p(k) \right) \geq 0, \,t\in\calT \label{eq.cond2wc}\nonumber\\
&\max_{[d]\in\calD}  B(t) \leq F(t), \;t\in\calT . \nonumber
\eea
%
%The statement of Proposition~\ref{pro:main} can be verified by evaluating explicitly
The worst-case quantities appearing in problem \eqref{eq:control1wc} can be evaluated explicitly as reported next; in the omitted derivations  we use repeatedly the fact that $\min_{|g|\leq r} h\tran g = -|h|\tran r$ and $\max_{|g|\leq r} h\tran g = |h|\tran r$:
 %\beas
 \[
 \begin{aligned}
 \underline p(t) &\doteq  \min_{|d(t)|\leq r(t)} p(t)   = \hat p(t) - |\Theta(t)|r(t) \\
 \underline u(t) &\doteq  \min_{|d(t)|\leq r(t)} u(t)   = \hat u(t) - |\Gamma(t)|r(t) \\
  \bar B(t)&\doteq  \max_{[d]\in\calD}  B(t) = \\
  & = \sum\nolimits_{\tau=0}^{t} \left(\one\tran \hat u(\tau) + \one\tran |\Gamma(\tau)|r(\tau)\right) \\
\bar J([\hat p],[\hat u], [\Theta], [\Gamma]) &\doteq
  \max_{[d]\in\calD} J([p],[ u]) \\
  & =
    J([\hat p],[\hat u])  +  \\
  & +\sum\nolimits_{t=0}^{T-1} \one\tran |\beta_t\Theta(t)  +\gamma\Gamma(t)  |r(t)  \\
  %&& + \sum_{t=0}^{T-1} \gamma  \one\tran |\Gamma(t)|r(t) \\
  \underline w(t+1) &\doteq \hat w(t+1) + \\
  & - \sum\nolimits_{k=0}^t |I + \Gamma(k) +(A\tran-I)\Theta(k)  |r(k) ,
  \end{aligned}
  %\eeas
  \]
 where for each $t\in\calT$ one has
 %\beas
  \[
 \begin{aligned}
\beta_t & \doteq  \eta \alpha^{T-t}+(1-\eta)a_t\\
 J([\hat p],[\hat u])  & \doteq \beta_0 \one\tran\bar p - \sum_{t=0}^{T-1}
 \beta_t\one\tran \hat p(t)    + \gamma \sum_{t=0}^{T-1}\one\tran \hat u(t), \\
 \hat w(t+1) & \doteq \sum_{k=0}^t \left( \hat e(k) +\hat u(k) + A\tran \hat p(k) - \hat p(k) \right).
 \end{aligned}
 \]
% \eeas
With the above positions, we can state the following

\begin{proposition}\label{pro:main}
The finite-horizon robust control problem \eqref{eq:control1wc} is equivalent to the   explicit linear program
\begin{align}
\min_{[\hat p],[\hat u], [\Theta], [\Gamma]} &\bar J([\hat p],[\hat u], [\Theta], [\Gamma])  \label{eq:control1wc2}\\
\mbox{s.t.:} &  [\underline p] \geq 0, \quad [\underline u]\geq 0 \nonumber \\
&\sum\nolimits_{k=0}^{t}  \alpha^{t-k}\left(\hat p(k) + |\Theta(k)|r(k)\right) \leq \alpha^t \bar p    ,\;   t\in\calT\nonumber\\
& \underline w(t+1)  \geq 0, \; t\in\calT  \nonumber\\
&\bar B(t) \leq F(t), \; t\in\calT . \nonumber
\end{align}
%\est
\end{proposition}

 \section{Numerical illustration}
 \label{sec:examples}
\begin{figure}[tb]
\centering
\includegraphics[width=.5\textwidth]{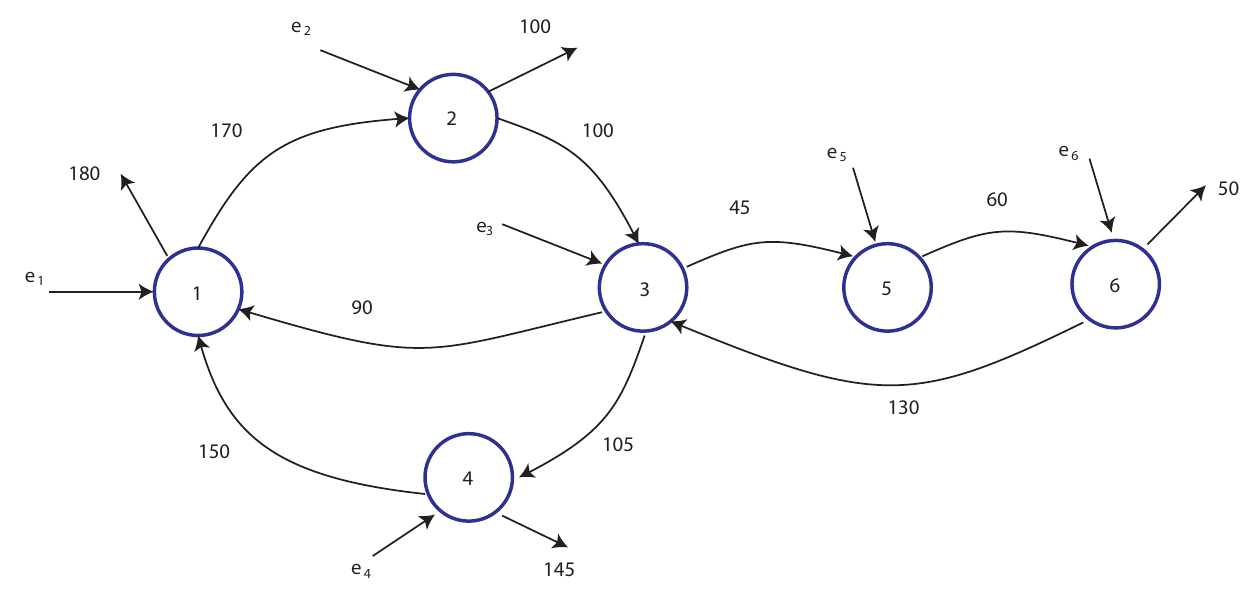}
\caption{A schematic network with $6$ nodes.}
\label{fig:6nodes}
\vspace{-0.5cm}
\end{figure}

 To illustrate the proposed approach, we consider a schematic network with $6$ nodes, plus the external fictitious node, as shown in
 Figure~\ref{fig:6nodes}.

The numbers on the edges in the graph in Figure~\ref{fig:6nodes} represent the initial nominal liabilities, forming the liability matrix $\bar P$,
vector $e=(e_1,\ldots,e_6)$ represents the external cash inflows at the nodes. We assume that the proportionality rule for default payments is in force.

\subsection{Nominal control}
Consider first a nominal scenario over a single period $T=1$, in which %a liquidity crisis makes all external in-flows drop by about 33\%, that
$\begin{smallmatrix} e(0)= (105,\, 25,\, 10,\, 190,\, 10,\, 120,\, 0) \end{smallmatrix}$. In this case, with no control,  all nodes default.
The clearing payments, computed according to \eqref{eq:clearingoptP_pr}, result to be
\[
P(0) =
% \left[\begin{array}{ccccccc}
\left[\begin{smallmatrix}0 & 164.6 & 0 & 0 & 0 & 0 & 174.3\\ 0 & 0 & 94.81 & 0 & 0 & 0 & 94.81\\ 86.18 & 0 & 0 & 100.5 & 43.09 & 0 & 0\\ 147.7 & 0 & 0 & 0 & 0 & 0 & 142.8\\ 0 & 0 & 0 & 0 & 0 & 53.09 & 0\\ 0 & 0 & 125.0 & 0 & 0 & 0 & 48.08\\ 0 & 0 & 0 & 0 & 0 & 0 & 0 %\end{array}\right].
\end{smallmatrix}\right].
\]
After such a clearing round, each node owes the residual amounts
$
(\begin{smallmatrix} 11.08   \;&10.38 \;&  10.18 \;&   4.45 \;&   6.91  \;&  6.91  \;& 0 \end{smallmatrix})
$,
for a total loss of $49.92$. The question now is the following: what could be a control intervention that may avoid the default?
To answer this question we solved the control problem \eqref{eq:control1}, over a single period $T=1$, setting parameters $\eta=0.9$, $\gamma =1$, and a total control budget $F(0) = 50$. The resulting optimal control action resulted to be $u(0) = (5,\;5,\;0,\; 0,\;5,\;0)$. It can be readily checked that with such control action the total in-flows are $c(0) = e(0) + u(0)$, and for such inputs the network returns to regular operations, with no default.
Overall, in this example, a relatively small intervention of amplitude $\|u(0)\|_1 = 15$ would be able to completely prevent the defaults and  bring the losses to zero.

We next consider a multi-stage setup with $T=3$ periods. We assume a $1\%$ interest rate on residual payments (i.e., $\alpha = 1.01$), and assume the following predicted stream of external payments
\[
e(0) = \left[\begin{smallmatrix}   105    \\ 0 \\   10 \\    0 \\    0 \\    0 \\    0  \end{smallmatrix}\right],
e(1) = \left[\begin{smallmatrix}    0  \\  25    \\ 0   \\190   \\  0  \\   0     \\0  \end{smallmatrix}\right],
e(2) = \left[\begin{smallmatrix}    0   \\  0    \\ 0   \\  0   \\ 10   \\120   \\  0  \end{smallmatrix}\right].
\]
We use, as in the previous case, $\eta=0.9$, $\gamma =1$, and we assume that the total control budget of $50$ is available
progressively as $F(0) = 15$, $F(1)=30$, $F(2) = 50$.
In this case, the solution of the control problem \eqref{eq:control1} gave us the optimal interventions
\[
u(0) = \left[\begin{smallmatrix}   2.19 \\    5.25 \\    0 \\    0\\    5.20 \\    2.36  \\         0 \end{smallmatrix}\right],
u(1) =
\left[\begin{smallmatrix}   2.84  \\   0  \\   0 \\    1.9  \\   0 \\  0  \\         0 \end{smallmatrix}\right],
u(2) = 0.
\]
These optimal injections, together with the computed optimal payment matrices at the intermediate times, are such that the network arrives to a regular (i.e., non default) situation at $T=3$. The optimal payment vectors were
\[
p(0) = \left[\begin{smallmatrix}
143.82 \\   75.11 \\ 61.32 \\   26.83 \\   16.69 \\   19.06 \\ 0
 \end{smallmatrix}\right],\;
p(1) =\left[\begin{smallmatrix}
131.04 \\   88.65 \\   51.27 \\  214.33 \\    9.61\\    9.61  \\ 0
 \end{smallmatrix}\right],\;
p(2) = \left[\begin{smallmatrix}
 77.96 \\   37.87 \\  130.49\\   57.09\\   34.47 \\  154.47     \\ 0
 \end{smallmatrix}\right].
\]
The full payment matrices can be deduced from the above payment vectors via the relation $P(t) = \diag(p(t))A$, where $A$ is the pro-rata matrix
\[
A = \left[\begin{smallmatrix}
 0 & 0.4857 & 0 & 0 & 0 & 0 & 0.5143\\ 0 & 0 & 0.5 & 0 & 0 & 0 & 0.5\\ 0.375 & 0 & 0 & 0.4375 & 0.1875 & 0 & 0\\ 0.5085 & 0 & 0 & 0 & 0 & 0 & 0.4915\\ 0 & 0 & 0 & 0 & 0 & 1.0 & 0\\ 0 & 0 & 0.7222 & 0 & 0 & 0 & 0.2778\\ 0 & 0 & 0 & 0 & 0 & 0 & 1
\end{smallmatrix}\right].
\]
The control effort in the present case amounts to a total $\one\tran(u(0)+u(1)+u(2)) = 19.74$, which is higher than the control effort needed in the single-stage case. This is expected since, due to interest, there is a price to pay for not having all the external payments available at $t=0$, and making the total control budget available only partially at the intermediate stages.

\subsection{Robust control}
We now examine the case of uncertain input flows.
We discuss first a single-step case ($T=1$). Consider the  nominal input cash flow
$\hat e(0) =(158,\,    38,\,    15,\,   285,\,    15,\,   180,\,     0)$.
In this nominal situation, the network  is in regular operation, all payments meet their liabilities,  no default occurs, and no corrective control action is needed.
 Assume, however, that the actual inputs are not exactly known, being however within a $33\%$ interval from the nominal values. By solving the robust control problem \eqref{eq:control1wc2} with  $\eta=0.9$, $\gamma =1$,
 and $F(0)=50$,
  we obtain
that  the optimal policies \eqref{eq:policy} are able to maintain the system default free in the worst case. This is achieved via the nominal control action and nominal payment
\[
\hat u(0) =  \left[\begin{smallmatrix}  1.32 \\    1.49\\    0.61 \\    0.13 \\   1.31 \\    0.65   \\ 0.0 \end{smallmatrix}\right],\quad
\hat p (0) =
\left[\begin{smallmatrix} 348.69 \\  198.52\\  239.40 \\  294.87 \\   58.70 \\  179.35 \\    0.0 \end{smallmatrix}\right],
\]
and reaction matrices
\beas
\Theta(0) &=&  \left[\begin{smallmatrix}
23 & 1.8 & 6.5 & 0.22 & 2.8 & 0.32 & 0\\ 5 & 95 & 1.1 & 0.049 & 0.46 & 0.054 & 0\\ 1.1 & 13 & 38 & 0.011 & 16 & 1.8 & 0\\ 0.11 & 0.75 & 1.9 & 0.93 & 1.7 & 0.16 & 0\\ 0.086 & 1 & 3.1 & 0 & 250 & 0.13 & 0\\ 0.029 & 0.34 & 1 & 0 & 55 & 6 & 0\\ 0 & 0 & 0 & 0 & 0 & 0 & 0
\end{smallmatrix}\right]\times 10^{-3}, \\
\Gamma(0) &=&
  \left[\begin{smallmatrix}
-23 & -1.8 & -6.6 & -0.22 & -2.8 & -0.32 & 0\\ -5.0 & -96 & -1.1 & -0.049 & -0.47 & -0.055 & 0\\ -1.1 & -13 & -39 & -0.011 & -16 & -1.8 & 0\\ -0.11 & -0.75 & -2 & -0.93 & -1.7 & -0.16 & 0\\ -0.087 & -1 & -3.1 & 0 & -250 & -0.13 & 0\\ -0.029 & -0.34 & -1 & 0 & -55 & -6.1 & 0\\ 0 & 0 & 0 & 0 & 0 & 0 & 0
\end{smallmatrix}\right]\times 10^{-3}.
\eeas
We finally consider a multi-step situation with $T=3$ and nominal external in-flows
\[
\hat e(0) = \left[\begin{smallmatrix}  15 \\ 30\\ 10\\ 200\\ 5\\ 100\\ 0   \end{smallmatrix}\right],
\hat e(1) = \left[\begin{smallmatrix}   80\\ 4\\ 0\\ 40\\ 5\\ 40\\ 0 \end{smallmatrix}\right],
\hat e(2) = \left[\begin{smallmatrix}   63\\ 4 \\ 5\\ 45\\ 5\\ 40\\ 0  \end{smallmatrix}\right].
\]
We assume that the external flow has 10\% uncertainty at $t=0$, while the ucertainty rises to 33\% at
$t=1$ and $t=2$. We let $\eta=0.9$, $\gamma =1$,  $\alpha = 1.01$, and $F(0) = 15$, $F(1)=30$, $F(2) = 50$.
Solving \eqref{eq:control1wc2} gives optimal policies that guarantee that the system is default free at the final time $T$, in all possible scenarios. The control effort was equal to $3.9$
in the nominal scenario and to $4.87$ in the  worst-case scenario,
 meaning that at most this sum is spent by the regulatory authority to maintain the system free of defaults.
 %in the worst case.

\section{Conclusions}
\label{sec:conc}
In this paper, we proposed a multi-period financial network model
 with an external control term representing corrective cash injections that can be performed by a ruling authority in order to prevent catastrophic cascaded failure events.
 %We considered both the general case where the payment matrices are unconstrained, and the pro-rata constrained case.
 We studied both the nominal
 case, in which  the cash inflows from the external sector are precisely known in advance, and  the more realistic case where the inflows are known only up to some interval of uncertainty. In this latter case, we proposed a robust approach based on linear feedback policies.
  In all the considered scenarios, the proposed control problems turn out to be efficiently solvable by means of linear programming.
   Numerical examples support the claim that small targeted interventions may avoid a cascaded failure effect and may thus significantly reduce the interbank contagion.

\bibliographystyle{IEEEtran}
\bibliography{biblio}
\newpage
\appendix
\section*{Proof of Lemma~\ref{lem.1}}

\subsection*{Statement 1}

To prove the first statement of Lemma, we first notice that if $([P^*],[u^*])$ is an optimal solution in~\eqref{eq:control1a}, then $[P^*]$ is an optimal solution in the problem
\bea
\min_{[P]} & J([P],[u^*])  \label{eq:control1a+}\\
\mbox{s.t.:} & [P] \geq 0,\nonumber \\
& \sum\nolimits_{k=0}^{t} \alpha^{t-k}P(k) \leq \alpha^t \bar P   ,\;   t\in\calT \nonumber\\
& C^*(t) + \sum\nolimits_{k=0}^{t}  \left(P(k) \tran- P(k) \right)\one \geq 0, \; t\in\calT\nonumber
\eea
(where $C^*(t)=c^*(0)+\ldots+c^*(t)$ corresponds to \emph{fixed}
control inputs $c^*(t)=e(t)+u^*(t)$, and the only decision variable is the sequence of payment matrices $[P]$).

Similarly, if $([p^*],[u^*])$ is an optimal solution in~\eqref{eq:control1}, then $[p^*]$ is an optimal solution in the problem
\bea
\min_{[p]} & J([p],[u^*])  \label{eq:control1+}\\
\mbox{s.t.:} & [p] \geq 0,\nonumber \\
& \sum\nolimits_{k=0}^{t} \alpha^{t-k}p(k) \leq \alpha^t \bar p   ,\;   t\in\calT \nonumber\\
& C^*(t) + \sum\nolimits_{k=0}^{t}  \left( A\tran p(k) - p(k) \right) \geq 0, \; t\in\calT\nonumber
\eea

Recalling the definition of cost functions~\eqref{eq:Ja} and~\eqref{eq:J}, one notices that the final budget $B(T-1)=B^*(T-1)$ is now also fixed, and hence
$J$ in~\eqref{eq:control1a+} (respectively,~\eqref{eq:control1+}) can be replaced by
\[
\tilde L([P]) = (1-\eta)L([P]) + \eta\one\tran \bar P(T) \one,
\]
or, respectively,
\[
\tilde L([p]) = (1-\eta)L([p]) + \eta\one\tran \bar p(T).
\]

The first statement of Lemma~\ref{lem.1} (absolute priority rule) is now implied\footnote{In the case of problem~\eqref{eq:control1a+}, our equation~\eqref{eq.no-defferal}
is a reformulation of~\cite[Equation~(30)]{AutPaper}, which is ensured by~\cite[Theorem~1]{AutPaper}. In the case of~\eqref{eq:control1a+}, the pro-rata rule entails that $p_{ji}^*(t)=a_{ji}p_j^*(t)$, so~\eqref{eq.no-defferal} is equivalent to~\cite[Equation~(30)]{AutPaper}, which is implied by~\cite[Theorem~1]{AutPaper}.}
by~\cite[Theorem~1]{AutPaper} (in the case of free payments) and~\cite[Theorem~1]{AutPaper} (in the case of pro-rata payments). Note that formally Theorems~1 and~2 in~\cite{AutPaper} are formulated for the special case $\eta=0$, where $\tilde L=L$ coincides with the total loss~\eqref{eq:Loss}. However, as noted in~\cite[Appendix~A.4]{AutPaper}, these theorems hold also for $\eta\in[0,1)$.\est

\subsection*{Statement 2}

Notice that~\eqref{eq.no-defferal}, in view of~\eqref{eq:wdynamics2free}, can be rewritten as follows: if $p_i(t)<\bar p_i^*(t)$ (some debt remains unpaid by period $t$), then $w_i^*(t)=0$. To prove the second statement, notice now that the banks whose liability has been paid by period $t$ (that is, $\bar p_i^*(t)=0$), obviously, do not receive additional cash at periods $t,t+1,\ldots,T-1$: otherwise, one could reduce the total budget $B(T-1)$ without violating any constraint.
Hence, if $u_i^*(t)>0$, then some debt remained unpaid ($p_i^*(s)<\bar p_i^*(s)$) at all periods $s=0,\ldots,t-1$, which, as has been noted, entails that $w_i^*(0)=\ldots=w_i^*(t)=0$. It remains to prove that $w_i^*(t+1)=0$, which will now be proved by contradiction.
Assume that $w_i^*(t+1)>0$. In view of~\eqref{eq.no-defferal}, one has
\[
w_i^*(t)+e_i^*(t)+u_i^*(t)+\sum_{j\ne i}p_{ji}^*(t)>\bar p_i^*(t).
\]
The latter inequality, however, remains valid if one reduces $u_i^*(t)$ by a small constant, decreasing thus also the total amount of case $B(T-1)$ and the value of cost function $J$. This contradicts to the solution's optimality.\est

\subsection*{Statement 3}

Note first that statement~3 follows from a formally weaker statement (A):

\emph{(A) For every optimal solution to the problem~\eqref{eq:control1a} or problem~\eqref{eq:control1} and every instant $t_0<T-1$ 
the implication holds: if $B^*(t_0)<F(t_0)$ (the budget constraint is not active at $t=t_0$), then $u^*(t_0+1)=0$.}

Indeed, suppose that $B^*(t_*)<F(t_*)$ yet $u^*(t)\ne 0$ at some instant $t\geq t_0+1$; let $t_1$ be the \emph{first} such instant.
Then, $t_1>t_*+1$ (due to statement (A)) and $u^*(t_*+1)=\ldots=u^*(t_1-1)=0$. Therefore, $B^*(t_1-1)=B^*(t_*)<F(t_*)\leq F(t_1-1)$. Statement (A) applied to
$t_0=t_1$ implies now that $u^*(t_1)=0$, which contradicts to the choice of $t_1$.

\subsubsection*{Proof of Statement (A)}

Assume that $B^*(t_0)<F(t_0)$ yet $u(t_0+1)\ne 0$. We will demonstrate that this assumption leads to the contradiction with the optimality of the solution, using the arguments similar to the ``advanced payment transformations'' from~\cite{AutPaper}.

Consider first a simpler case of free payments (the problem~\eqref{eq:control1}). Let $i$ be one of the banks that receive extra cash at time $t_0+1$: $u_i(t_0+1)>0$. Due to statement~2, this is possible only when $\bar p_i^*(t_0+1)>0$ and, furthermore,~\eqref{eq.no-defferal} 
implies (in view of $c_i^*(t_0+1)\geq u_i^*(t_0+1)>0$) that $p_i^*(t_0+1)>0$, so at least bank $j\ne i$ receives payment from $i$: 
$p_{ij}^*(t_0+1)>0$, and hence $p_{ij}^*(t_0)<\bar p_{ij}^*(t_0)$.

Define the sequence $[P]$ of payment matrices as follows
\[
p_{km}(t)=
\begin{cases}
p_{ij}^*(t_0)+\alpha^{-1}\varepsilon,& (i,j)=(k,m),t=t_0,\\
p_{ij}^*(t_0+1)-\varepsilon,& (i,j)=(k,m),t=t_0+1,\\
p_{km}^*(t),&\text{in all other cases}.
\end{cases}
\]
and also a new sequence of control inputs $[u]$, where
\[
u_k(t)=
\begin{cases}
u_i^*(t_0)+\alpha^{-1}\varepsilon, & k=i, t=t_0,\\
u_i^*(t_0+1)-\varepsilon, & k=i, t=t_0+1,\\
u^*_j(t_0+1)+(1-\alpha^{-1})\varepsilon,& k=j,\,t=t_0+1\\
u^*_k(t),& \text{in all other cases}.
\end{cases}
\]
In other words, $i$ pays to $j$ a larger amount at time $t=t_0$ in order to decrease the payment at time $t=t_0+1$. To make this possible without violating the inequality $w(t)\geq 0$, one has to increase the cash input injected to $i$ at time $t=t_0$ (to make  and the cash input to $j$ at time $t=t_0+1$; at the same time, the cash input to $i$ at time $t=t_0$ have to be decreased in order to preserve the total budget.
Here $\varepsilon>0$ is such that
\[
\begin{gathered}
p_{ij}(t_0+1)>0, u_i(t_0+1)>0, p_{ij}(t_0)<\bar p_{ij}^*(t_0),\\
 B(t_0)=B^*(t_0)+\alpha^{-1}\ve<F(t_0).
\end{gathered}
\]
By construction, $P(t),u(t)$ are nonnegative. One may also notice that $([P],[u])$ satisfy the conditions~\eqref{eq.cond1a},~\eqref{eq.cond2a} and~\eqref{eq.cond3a}.
The condition~\eqref{eq.cond3a} at time $t=t_0$ is guaranteed by the choice of $\varepsilon$, whereas
$B(t)=B^*(t)\leq F(t)\,\forall t\ne t_0$. 

The condition~\eqref{eq.cond1a} is not violated for $t\ne t_0$, because
the its left-hand side remains invariant after the replacement of $[P^*]$ by $[P]$. To verify this condition at $t=t_0$,
recall that~\eqref{eq.cond1a} is nothing else than the inequality $P(t)\leq\bar P(t)$, which holds at $t=t_0$ by construction.

The condition~\eqref{eq.cond2a} is equivalent to the relation $w_k(t+1)\geq 0\,\forall k=1,\ldots,n$. It can be easily seen that, by construction, one has $w_k(t+1)=w_k^*(t+1)$ for all $k\ne j$, whereas
\[
w_j(t+1)=
\begin{cases}
w_j^*(t+1),& t\ne t_0\\
w_j^*(t_0+1)+\alpha^{-1}\varepsilon>w_j^*(t_0+1),& t=t_0.
\end{cases}
\]
Hence,~\eqref{eq.cond2a} also holds. At the same time,
\[
J([P],[u])=J([P^*],[u^*])+(\alpha^{-1}a_{t_0}-a_{t_0+1})\ve<J([P^*],[u^*]),
\]
which leads us to the contradiction with the optimality of $([P^*],[u^*])$. Hence, the statement (A) is valid.

The case of pro-rata constraint is considered similarly with the only difference that, transferring the payment of bank $i$ from period $t_0+1$ to $t_0$,
the control intervention at time $t=t_0+1$ is needed by all banks $j\ne i$. Instead of sequence of matrices $[P]$, one can 
construct a sequences of payment vectors $[p]$ by defining
\[
p_{k}(t)=
\begin{cases}
p_{i}^*(t_0)+\alpha^{-1}\varepsilon,& k=i\;\text{and}\;t=t_0,\\
p_{i}^*(t_0+1)-\varepsilon,& k=i\;\text{and}\;t=t_0+1,\\
p_{k}^*(t),&\text{in all other cases}.
\end{cases}
\]
and a sequence of control inputs
\[
u_k(t)=
\begin{cases}
u_i^*(t_0)+\alpha^{-1}\varepsilon, & k=i, t=t_0,\\
u_i^*(t_0+1)-\varepsilon, & k=i, t=t_0+1,\\
u^*_k(t_0+1)+(1-\alpha^{-1})\varepsilon a_{ik},& k\ne i,\,t=t_0+1\\
u^*_k(t),& \text{in all other cases}.
\end{cases}
\]
Here $\ve>0$ is so small that $p_i(t_0)<\bar p_i^*(t_0)$, $p_i(t_0+1)>0$, $u_i^*(t_0+1)>0$ and $B(t_0)=B^*(t_0)+\alpha^{-1}\ve<F(t_0)$. It can be shown that $J([p],[u])<J([p^*],[u^*])$ and the pair of sequences $([p],[u])$
is feasible, in particular, $w_k(t+1)=w_k^*(t+1)$ for all $t\ne t_0$ and all $k$ and
\[
w_{k}(t_0+1)=
\begin{cases}
w_i^*(t_0+1),& k=i\\
w_k^*(t_0+1)+\alpha^{-1}\varepsilon a_{jk}>w_k^*(t_0+1),& k\ne i.
\end{cases}
\]
This leads to the contradiction with optimality of $([p^*],[u^*])$. Statement (A) is proved.\est
\end{document}